\documentclass{article}
\usepackage{graphicx} 
\usepackage[utf8]{inputenc}

\usepackage{amsthm}
\usepackage{amsmath}
\usepackage{amssymb}
\usepackage{url}
\usepackage{float}
\usepackage{bm}
\usepackage{xcolor}
\usepackage{graphicx}
\usepackage{mathtools}
\usepackage{subfigure}

\graphicspath{{../}}

\newcommand{\HH}{\mathcal{H}}
\newcommand{\scal}[2]{{\left\langle{{#1}\mid{#2}}\right\rangle}}
 
\newcommand{\Id}{\ensuremath{\operatorname{Id}}}

\newcommand{\prox}{{\rm prox}}

\newcommand{\RR}{\ensuremath{\mathbb{R}}}
\newcommand{\NN}{\ensuremath{\mathbb{N}}}

\newcommand{\RX}{\ensuremath{\left]-\infty,+\infty\right]}}

\newtheorem{problem}{Problem}[section]

\newtheorem{remark}{Remark}[section]
\newtheorem{theorem}{Theorem}[section]
\newtheorem{corollary}{Corollary}[section]
\newtheorem{proposition}{Proposition}[section]

\title{Lyapunov analysis for FISTA under strong convexity}
\author{Luis M. Briceño-Arias\footnote{Universidad Técnica 
Federico Santa María, Department of Mathematics, 
luis.briceno@usm.cl. This work was supported by the National 
Agence of Research and Development (ANID) from Chile, under 
the grants FONDECYT 1230257, MATH-AmSud 23-MATH-17, 
{MATH-AmSud 250023,} and 
Centro de Modelamiento Matemático (CMM) BASAL fund 
FB210005 for centers of excellence.}}
\date{}

\begin{document}
\maketitle

\begin{center}
\textit{To Hédy, in deep gratitude for his inspiring contributions to 
mathematics and the Chilean optimization community.}
\end{center}
\begin{abstract}
In this paper, we conduct a theoretical and numerical study of the 
Fast Iterative Shrinkage-Thresholding Algorithm (FISTA) 
under 
strong convexity assumptions. We propose an autonomous 
Lyapunov function that reflects the strong convexity of the objective 
function, whether it arises from the smooth or non-smooth 
component. This Lyapunov function decreases monotonically at a 
linear rate along the iterations of the algorithm for a fixed inertial 
parameter{, provided that the inertial parameter is built from 
the total strong convexity modulus of the objective function}. 
{Our analysis shows that the resulting algorithm with optimal 
parameters, its Lyapunov function, and its linear convergence rate 
are invariant under any transfer of strong convexity between the 
smooth and non-smooth components. This invariance also holds 
for forward-backward splitting (FBS).}
Within this framework, we compare the performance of 
{FBS} and some FISTA
variants, and 
find that this strategy leads FISTA to outperform all other 
configurations, including FBS. Moreover, we identify parameter 
regimes in which FBS yields better performance than FISTA 
when 
the strong convexity of the non-smooth part is not leveraged 
appropriately.
\end{abstract}

\section{Introduction}
In this paper, we provide a theoretical and numerical comparison 
of inertial and non-inertial algorithms for solving the optimization 
problem
\begin{equation}
\label{e:base}
\min_{x\in\HH}f(x)+h(x),
\end{equation}
where $\HH$ is a real Hilbert space, $f\colon\HH\to\RR$
is a convex differentiable function such that $\nabla f$ is 
$L$-Lipschitz continuous, $h\colon\HH\to\RX$ is proper, lower 
semicontinuous, and convex, and we assume the existence of 
solutions. This problem is pervasive in signal and image 
processing problems 
\cite{BTeb,Sigpro1,ChPockrev,combettes2005signal,%
daubechies2004iterative},
 mean field games
\cite{benamou2015augmented, svva}, among other 
disciplines (see, e.g., \cite{Livre1,combettes2005signal} and the 
references therein). 

A classical algorithm for solving \eqref{e:base} is 
the forward-backward splitting (or 
proximal-gradient algorithm), which iterates
\begin{equation}
\label{e:FBS}\tag{FBS}
x_0\in\HH,\quad (\forall k\in\NN)\quad x_{k+1}={\rm 
prox}_{\gamma 
h}(x_k-\gamma\nabla f(x_k)).
\end{equation}
The convergence to a solution to \eqref{e:base} of the sequence 
$(x_k)_{k\in\NN}$ generated by \eqref{e:FBS} is guaranteed 
provided that 
$\gamma\in\left]0,2/L\right[$. This algorithm finds its roots in the 
projected gradient method \cite{levitin1966constrained}
 and certain methods for solving variational inequalities 
\cite{auslender1972problemes, brezis1967methodes, 
mercier1979topics, sibony1970methodes}.
{In image processing, }
 when $h$ is the $\ell^1$ norm and $f$ is a quadratic data fidelity 
term, \eqref{e:FBS} is widely known as ISTA 
\cite{daubechies2004iterative}{, whose proximal-gradient 
structure 
is derived in \cite{combettes2005signal}.}
In the case when $f$ is $\mu$-strongly convex for some $\mu>0$, 
\eqref{e:FBS} achieves 
a linear convergence rate \cite{TayFBS}, which is 
optimized for $\gamma^*=2/(L+\mu)$. This optimal step-size for 
the strongly convex case helps to understand the standard choice 
of $\gamma\approx 2/L$ in the absence of strong convexity 
\cite{Sigpro1}.

In the case when $h=0$, \eqref{e:FBS} reduces to the classical 
gradient descent method, which 
was accelerated from $O(1/k)$ {(in fact $o(1/k)$, see
\cite{CombettesSalzoVilla2018})} to 
$O(1/k^2)$ in the convergence in values by considering an 
additional inertial step with a specific choice of inertial parameters 
$(\alpha_k)_{k\in\NN}\subset{\left]0,1\right[}$ \cite{Nest83} 
\begin{equation}
\label{e:nest}
(\forall k\in\NN)\quad 
\begin{cases}
x_{k+1}=y_k-\gamma\nabla f(y_k)\\
y_{k+1}=x_{k+1}+\alpha_k(x_{k+1}-x_k),
\end{cases}
\end{equation}
where $x_0$ and $y_0$ are in $\HH$. This accelerated method 
was studied 
as a discretization of a second order dynamical system in 
\cite{suBC}, where the authors also remarked the non-monotone 
behavior of errors in function values and iterates (see also 
\cite{apid}). 

The \textit{fast iterative soft-thresholder algorithm} (FISTA) 
algorithm \cite{BTeb} extends this accelerated 
scheme to the case $h\neq 0$ inspired from \cite{gul92,Nest83} 
via the recurrence
\begin{equation}
\label{e:fista}\tag{FISTA}
(\forall k\in\NN)\quad 
\begin{cases}
x_{k+1}={\rm prox}_{\gamma h}(y_k-\gamma\nabla f(y_k))\\
y_{k+1}=x_{k+1}+\alpha_k(x_{k+1}-x_k),
\end{cases}
\end{equation}
where $x_0$ and $y_0$ are arbitrary in $\mathcal{H}$, 
$(\alpha_k)_{k\in\NN}$ is a sequence of inertial 
parameters {in $\left]0,1\right[$}, and $\gamma=1/L$. In 
\cite{BTeb}, the convergence 
of 
values is guaranteed but the convergence of the iterates remained 
an open question.
In \cite{attchb,AttP16,ChDoss} the convergence of the iterates 
is guaranteed for a suitable choice of parameters  
$(\alpha_k)_{k\in\NN}\subset{\left]0,1\right[}$ {(see also
\cite{CombettesGlaudin2017} for a study of more general inertial 
parameters)}. In particular, 
the convergence of the iterates is demonstrated in
\cite{attchb} via a Lyapunov analysis of a 
continuous dynamical system inspired
from \cite{AttAl01}. Similar techniques are used in
\cite{AttP16} to obtain the better convergence rate of $o(1/k^2)$.
However, as perceived in \cite{attouch2016fast}, 
\cite[Chapter~10.7.4]{Beck}, and 
\cite{suBC}, values and iterates of the sequences generated by 
\eqref{e:fista} exhibit an oscillatory asymptotic behavior, given that 
it is not a descent method.

In the particular case when $f$ is $\mu$-strongly convex for some 
$\mu>0$, similar techniques to those used in dynamical systems 
are used by studying particular energy Lyapunov functionals  
exploiting the geometry 
of the problem to achieve the best convergence rate.
In \cite[Theorem~4.10]{ChPockrev}, 
\cite[Theorem~6]{aujol2024fista}, 
\cite[Corollary~4.21]{daspremont2021acceleration}, 
\cite[Theorem~9]{suBC}
non-autonomous Lyapunov 
functions for the study of \eqref{e:fista} with adaptative inertia are 
proposed. In \cite[Theorem~9]{suBC} a polynomial rate of 
$O(1/k^3)$ is obtained, while linear convergence rates are derived 
in \cite{aujol2024fista,ChPockrev,daspremont2021acceleration}.
In particular, the result 
in \cite[Theorem~4.10]{ChPockrev} takes advantage of the case 
when $h$ is $\rho$-strongly convex for some $\rho\ge 0$ while 
a quadratic growth condition is assumed in 
\cite[Theorem~6]{aujol2024fista}.
On the other hand, a variant of \eqref{e:fista} with fixed inertia 
parameter 
known as V\ref{e:fista} leads to a linear convergence rate of 
$O((1-\sqrt{\mu/L})^k)$ (see, e.g.,  \cite[Theorem~10.42]{Beck}).
In \cite[Theorem~4.14]{daspremont2021acceleration} a similar 
rate is obtained in the case $h=0$. These approaches do not 
leverage the strong convexity of $h$.

In this paper, we generalize the best-known theoretical linear 
convergence rates in the literature when both $f$ and $h$ are 
strongly convex for \eqref{e:fista} with a {constant} inertia 
parameter via 
autonomous Lyapunov functions. Even if these assumptions may 
be considered restrictive, they provide valuable theoretical insights 
into the structure of the algorithm and the role of inertia in its 
behavior. {We prove that \eqref{e:fista} achieves a linear rate 
for the constant inertial parameter}
$${\alpha=\frac{\sqrt{L+\rho}-\sqrt{\mu+\rho}}
{\sqrt{L+\rho}+\sqrt{\mu+\rho}}}$$
{leading to a linear convergence rate of}
$$r=1-\sqrt{\frac{\mu+\rho}{L+\rho}}.$$
This convergence rate 
coincides with the best-known linear convergence rate in 
\cite{ChPockrev}, but with the difference that we use a constant 
inertial parameter and an 
autonomous Lyapunov function. 
{When $\rho>0$, this rate is strictly lower than 
$1-\sqrt{\mu/L}$, which is the 
optimal rate of \eqref{e:fista} when the strong convexity of the 
non-smooth part is not
exploited. We also show that the algorithm with optimal 
parameters, its 
Lyapunov function, and its rate are invariant under any transfer of 
strong convexity between $f$ and $h$. This invariance also holds 
for forward-backward splitting (FBS). }
We also provide a 
Lyapunov-type analysis for the linear convergence of 
\eqref{e:fista}, which allows us to 
study a decreasing quantity explaining the asymptotic behavior of 
inertial algorithms without oscillations. From this Lyapunov 
analysis, we derive complementary convergence rates for iterates 
and values and we illustrate our theoretical results via numerical 
experiments in an academic example, comparing inertial and 
non-inertial approaches.

The paper is organized as follows. Section~\ref{sec:notation} 
presents a concise overview of the notation and preliminary 
concepts used throughout the paper. In Section~\ref{sec:main}, we 
present our main result, introducing the proposed autonomous 
Lyapunov function along with the asymptotic analysis that 
establishes its linear convergence rate. Section~\ref{sec:lev} 
explores {the effect of transferring} the strong convexity 
{between} the functions $f$ and $h$ in both \eqref{e:FBS} and 
\eqref{e:fista}. Finally, 
Section~\ref{sec:num} provides 
numerical experiments that illustrate the theoretical and numerical 
bounds derived in the previous sections, using an academic 
example.

\section{Notation}
\label{sec:notation}
Throughout this paper, $\HH$ is a real Hilbert space, endowed 
with scalar product $\scal{\cdot}{\cdot}$ and induced norm 
$\|\cdot\|=\sqrt{\scal{\cdot}{\cdot}}${, and $\Id$ denotes the 
identity operator on $\HH$}.
Let $f\colon\HH\to\RX$ be a function. Given $\rho\ge 0$, $f$
is $\rho$-strongly convex if $f-\rho\|\cdot\|^2/2$ is convex. In 
particular, convex functions are $0$-strongly convex. We denote 
by $\Gamma_0(\HH)$ the class of proper lower semicontinuous 
convex functions. The problem of minimizing a 
function $f\in\Gamma_0(\HH)$ which is $\rho$-strongly convex, for 
some $\rho>0$, over the space $\HH$ admits a unique solution.
Suppose that $f\in\Gamma_0(\HH)$. The proximity operator
of $f$ is
\begin{equation}
\prox_{f}\colon\HH\to\HH\colon 
x\mapsto\arg\min_{y\in\HH}f(y)+\frac{1}{2}\|y-x\|^2,
\end{equation}
which is well defined since the objective function above is 
$1$-strongly convex, proper, and lower semicontinuous.
For every $\gamma >0$, $x$ and $p$ in $\HH$, we have
\begin{equation}
\label{e:caractprox}
p=\prox_{\gamma f}x\quad\Leftrightarrow\quad 
\frac{x-p}{\gamma}\in\partial f(p),
\end{equation}
where $\partial f$ denotes the subdifferential of $f$. In the case 
when $f$ is Fréchet-differentiable at some $x\in\HH$, $\partial 
f(x)=\{\nabla f(x)\}$, where 
$\nabla f$ is the 
gradient operator of $f$. For further results in convex analysis, the 
reader is referred to \cite{Livre1, nesterov2004introductory}.

\section{Main results}
\label{sec:main}
Consider the following optimization problem.
\begin{problem}
\label{prob:1}
Let $L>0$, $\mu\in\left[0,L\right[$, and $\rho\ge 0$ be such that 
$\mu+\rho>0$. Let 
$f\colon\mathcal{H}\to\mathbb{R}$
be a differentiable $\mu$-strongly convex function such that 
$\nabla f$
is $L$-Lipschitz continuous.
Let $h\in\Gamma_0(\mathcal{H})$ be $\rho$-strongly convex.
The problem is to
\begin{equation}
\label{e:prob1}
\min_{x\in\mathcal{H}}F(x):=f(x)+h(x),
\end{equation}
which has a unique solution $x^*$.
\end{problem}
 In this section
we provide a Lyapunov function adapted to the hypotheses in 
Problem~\ref{prob:1} to study the linear convergence rate of 
\eqref{e:fista} {with a constant inertial parameter 
$\alpha_k\equiv\alpha\in\left]0,1\right[$}. First, given 
$c\in\left]0,1\right[$, $\widetilde{x}_0$ 
and 
$\widetilde{z}_0$ in $\HH$, 
\eqref{e:fista} {with constant inertia} can be equivalently 
written as 
\begin{equation}
\label{e:fista_z}
(\forall k\in\NN)\quad 
\begin{cases}
\widetilde{y}_k=\widetilde{x}_k+c(\widetilde{z}_k-\widetilde{x}_k)\\
\widetilde{x}_{k+1}={\rm prox}_{\gamma h}
(\widetilde{y}_k-\gamma\nabla f(\widetilde{y}_k))\\
\widetilde{z}_{k+1}=\frac{\alpha}{1-c}\widetilde{z}_k
-\frac{\alpha}{c(1-c)}\widetilde{y}_k
+\frac{c+\alpha}{c}\widetilde{x}_{k+1}
\end{cases}
\end{equation}
as the following result asserts. {The parameter 
$c\in\left]0,1\right[$ is free at this stage and will be determined 
below in 
order to optimize the linear convergence rate of \eqref{e:fista}.}
\begin{proposition}
\label{p:1}
{Let $c\in\left]0,1\right[$ and let $\alpha\in\left]0,1\right[$.}
Let $(x_k)_{k\in\NN}$ and $(y_k)_{k\in\NN}$
be the sequences generated by the algorithm in \eqref{e:fista} 
{with $\alpha_k\equiv\alpha$},
for some $x_0$ and $y_0$ in $\HH$
and let $(\widetilde{x}_k)_{k\in\NN}$, 
$(\widetilde{y}_k)_{k\in\NN}$, and $(\widetilde{z}_k)_{k\in\NN}$  
be 
the sequences generated by the algorithm in \eqref{e:fista_z} for 
some $\widetilde{x}_0$ and $\widetilde{z}_0$ in $\HH$.
Suppose that 
$${y}_0=\widetilde{x}_0+c(\widetilde{z}_0-\widetilde{x}_0)
\quad\text{and}\quad 
{x}_0=\widetilde{x}_0.$$
Then, {by setting}
\begin{equation}
\label{e:defzk}
{(\forall k\in\NN)\quad 
z_k:=x_k+\frac{1}{c}(y_k-x_k),}
\end{equation}
{we have, for every $k\in\NN$,} $x_k=\widetilde{x}_k$, 
$y_k=\widetilde{y}_k$, and
${z_k=}\widetilde{z}_k$.
\end{proposition}
\begin{proof}
We will proceed by recurrence. 
Note that, by definition, ${y}_0=\widetilde{y}_0$, 
${x}_0=\widetilde{x}_0$ and
$\widetilde{z}_0=\widetilde{x}_0+\frac{1}{c}(y_0-\widetilde{x}_0)=
{x}_0+\frac{1}{c}(y_0-{x}_0){=z_0}$.
Now set $k\in\NN$, and suppose that, for every $\ell\le k$, 
$x_{\ell}=\widetilde{x}_{\ell}$,
$y_{\ell}=\widetilde{y}_{\ell}$, and 
$\widetilde{z}_\ell={x}_\ell+\frac{1}{c}(y_\ell-x_\ell)
{=z_\ell}$. 
Then, it follows from \eqref{e:fista_z} and 
\eqref{e:fista} that $x_{k+1}=\widetilde{x}_{k+1}$ and
\begin{align}
\widetilde{y}_{k+1}&=\widetilde{x}_{k+1}+c(\widetilde{z}_{k+1}-\widetilde{x}_{k+1})\nonumber\\
&=(1+\alpha)\widetilde{x}_{k+1}+\frac{\alpha}{1-c}(c\widetilde{z}_k
-\widetilde{y}_k)\nonumber\\    
&=(1+\alpha)\widetilde{x}_{k+1}-\alpha \widetilde{x}_k\\
&=(1+\alpha){x}_{k+1}-\alpha {x}_k={y}_{k+1}.
\end{align}
Hence, it follows from \eqref{e:fista_z} that $\widetilde{z}_{k+1}=
\widetilde{x}_{k+1}+\frac{1}{c}(\widetilde{y}_{k+1}-\widetilde{x}_{k+1})
={x}_{k+1}+\frac{1}{c}({y}_{k+1}-{x}_{k+1}){=z_{k+1}}$, 
which completes the 
proof.
\end{proof}

For every $\mu\ge 0$ {and $\rho\ge 0$},  we define 
\begin{equation}
\label{e:deflyap}
{\Phi_{\mu,\rho}}\colon (x,z)\mapsto 
F(x)-F(x^*)+{\frac{\mu+\rho}{2}}
\|z-x^*\|^2.
\end{equation}
The following is our main result, which asserts that 
$({\Phi_{\mu,\rho}}(x_k,z_k))_{k\in\NN}$ satisfies a 
Lyapunov-type 
property as in \cite{AttP16}.  
\begin{theorem}
\label{t:1}
In the context of Problem~\ref{prob:1},
let $(x_k)_{k\in\NN}$ and $(y_k)_{k\in\NN}$ be the sequences 
generated by the algorithm \eqref{e:fista}, for some $x_0$ and 
$y_0$ in $\HH$
and {constant} parameters 
\begin{equation}
\label{e:param}
 \gamma=\frac{1}{L}\quad\text{and}\quad 
{\alpha=\frac{\sqrt{L+\rho}-\sqrt{\mu+\rho}}
{\sqrt{L+\rho}+\sqrt{\mu+\rho}}}.
\end{equation}
Then, by setting
\begin{equation}
\label{e:defz}
(\forall k\in\NN)\quad 
z_k={x}_k+{\frac{\sqrt{L+\rho}+\sqrt{\mu+\rho}}
{\sqrt{\mu+\rho}}}
(y_k-x_k),
\end{equation}
we have
\begin{equation}
\label{e:goal0}
(\forall k\in\NN)\quad {\Phi_{\mu,\rho}}(x_{k+1},z_{k+1})\le 
r(\mu,\rho,L){\Phi_{\mu,\rho}}(x_{k},z_{k}),
\end{equation}
where 
\begin{equation}
\label{e:defr}
r(\mu,\rho,L)={1-\sqrt{\frac{\mu+\rho}{L+\rho}}}.
\end{equation}
\end{theorem}
\begin{proof}
Fix $k\in\NN$. It follows from \eqref{e:caractprox} and 
\eqref{e:fista} that
\begin{equation}
\frac{y_k-x_{k+1}}{\gamma}-\nabla f(y_k)\in\partial h(x_{k+1}).
\end{equation}
Then, the $\rho$-strong convexity 
of $h$ implies
\begin{multline}
\label{e:hstrong1}
0\ge h(x_{k+1})-h(x^*)+\frac{1}{\gamma}\scal{y_k-x_{k+1}}{x^*-x_{k+1}}-\scal{\nabla f(y_k)}{x^*-x_{k+1}}\\
+\frac{\rho}{2}\|x_{k+1}-x^*\|^2
\end{multline}
 and 
 \begin{multline}
 \label{e:hstrong2}
0\ge h(x_{k+1})-h(x_k)+\frac{1}{\gamma}\scal{y_k-x_{k+1}}{x_k-x_{k+1}}-\scal{\nabla f(y_k)}{x_k-x_{k+1}}\\
+\frac{\rho}{2}\|x_{k+1}-x_k\|^2.
\end{multline}
Moreover, the $\mu$-strong convexity of $f$ yields
 \begin{equation}
 \label{e:fstrong1}
0\ge f(y_{k})-f(x^*)+\scal{\nabla f(y_k)}{x^*-y_k}+\frac{\mu}{2}\|y_{k}-x^*\|^2
\end{equation}
and 
\begin{equation}
\label{e:fstrong2}
0\ge f(y_{k})-f(x_k)+\scal{\nabla f(y_k)}{x_k-y_k}+\frac{\mu}{2}\|y_{k}-x_k\|^2.
\end{equation}
By adding \eqref{e:hstrong1} and \eqref{e:fstrong1} we obtain
\begin{multline}
\label{e:comb1}
0\ge h(x_{k+1})+f(y_k)-F(x^*)+\frac{1}{\gamma}\scal{y_k-x_{k+1}}{x^*-x_{k+1}}+\scal{\nabla f(y_k)}{x_{k+1}-y_k}\\
+
\frac{\rho}{2}\|x_{k+1}-x^*\|^2+\frac{\mu}{2}\|y_{k}-x^*\|^2
\end{multline}
and, analogously, by adding \eqref{e:hstrong2} and \eqref{e:fstrong2} we get
\begin{multline}
\label{e:comb2}
0\ge h(x_{k+1})+f(y_k)-F(x_k)+\frac{1}{\gamma}\scal{y_k-x_{k+1}}{x_k-x_{k+1}}+\scal{\nabla f(y_k)}{x_{k+1}-y_k}\\
+
\frac{\rho}{2}\|x_{k+1}-x_k\|^2+\frac{\mu}{2}\|y_{k}-x_k\|^2.
\end{multline}
On the other hand, it follows from the $L$-Lipschitz continuity of $\nabla f$ that
\begin{equation}
\label{e:Llip}
0\ge f(x_{k+1})-f(y_k)-\scal{\nabla f(y_k)}{x_{k+1}-y_k}
-\frac{L}{2}\|x_{k+1}-y_k\|^2.
\end{equation}
Now let $\lambda_1$, $\lambda_2$, and $\lambda_3$ be strictly positive reals such that $\lambda_3=\lambda_1+\lambda_2$. Then by multiplying 
\eqref{e:comb1} by $\lambda_1$, \eqref{e:comb2} by $\lambda_2$, and \eqref{e:Llip} by $\lambda_3$ and adding
the three resulting inequalities, we obtain
\begin{align}
\label{e:super}
0&\ge \lambda_3 (F(x_{k+1})-F(x^*))-\lambda_2 (F(x_{k})-F(x^*))\nonumber\\
&\qquad+\frac{\lambda_1}{\gamma}\scal{y_k-x_{k+1}}{x^*-x_{k+1}}+\frac{\lambda_2}{\gamma}\scal{y_k-x_{k+1}}{x_k-x_{k+1}}\nonumber\\
&\qquad+\frac{\lambda_1\rho}{2}\|x_{k+1}-x^*\|^2+\frac{\lambda_1\mu}{2}\|y_{k}-x^*\|^2\nonumber\\
&\qquad+\frac{\lambda_2\rho}{2}\|x_{k+1}-x_k\|^2+\frac{\lambda_2\mu}{2}\|y_{k}-x_k\|^2-\frac{\lambda_3L}{2}\|x_{k+1}-y_k\|^2\nonumber\\
&=\lambda_3 (F(x_{k+1})-F(x^*))-\lambda_2 (F(x_{k})-F(x^*))\nonumber\\
&\qquad-\frac{\lambda_1}{2}\left(\frac{1}{\gamma}-\mu\right)\|y_{k}-x^*\|^2-\frac{\lambda_2}{2}\left(\frac{1}{\gamma}-\mu\right)\|y_{k}-x_k\|^2\nonumber\\
&\qquad +\frac{\lambda_1}{2}\left(\rho+\frac{1}{\gamma}\right)\|x_{k+1}-x^*\|^2+\frac{\lambda_2}{2}\left(\rho+\frac{1}{\gamma}\right)\|x_{k+1}-x_k\|^2\nonumber\\
&\qquad +\frac{\lambda_3}{2}\left(\frac{1}{\gamma}-L\right)\|x_{k+1}-y_k\|^2\nonumber\\
&=\lambda_3 (F(x_{k+1})-F(x^*))-\lambda_2 (F(x_{k})-F(x^*))\nonumber\\
&\qquad -\frac{\lambda_1}{2}\left(\frac{1}{\gamma}-\mu\right)\|y_{k}-x^*\|^2-\frac{\lambda_2c^2}{2}\left(\frac{1}{\gamma}-\mu\right)\|z_{k}-x_k\|^2\nonumber\\
&\qquad +\frac{\lambda_1}{2}\left(\rho+\frac{1}{\gamma}\right)\|x_{k+1}-x^*\|^2+\frac{\lambda_2}{2}\left(\rho+\frac{1}{\gamma}\right)\|x_{k+1}-x_k\|^2\nonumber\\
&\qquad 
+\frac{\lambda_3}{2}\left(\frac{1}{\gamma}-L\right)\|x_{k+1}-y_k\|^2,
\end{align}
where we have used {that $y_k-x_k=c(z_k-x_k)$, which 
follows from the definition \eqref{e:defzk} of $z_k$ in 
Proposition~\ref{p:1}, for an arbitrary value of $c\in\left]0,1\right[$ 
to be fixed below.}
Note that we may assume, without loss of generality, that 
$\lambda_2=1$. Indeed, otherwise we can divide the last 
inequality by $\lambda_2$ and consider the new variables 
$\widetilde{\lambda}_1=\lambda_1/\lambda_2$ and 
$\widetilde{\lambda}_3=\lambda_3/\lambda_2$.
{On the other hand, by recalling the two elementary 
identities}
\begin{equation}
\label{e:identz}
{z_k-y_k=(1-c)(z_k-x_k)}
\end{equation}
{and}
\begin{equation}
\label{e:identx}
{x_{k+1}-y_k=(1-c)(x_{k+1}-x_k)+c(x_{k+1}-z_k),}
\end{equation}
{which follow directly from $y_k=x_k+c(z_k-x_k)$,  
it follows from \eqref{e:fista_z}
and \cite[Lemma~2.14]{Livre1} that}
\begin{align*}
\|z_{k+1}-x^*\|^2&=\Big\|\frac{\alpha}{1-c}(z_k-x^*)-\frac{\alpha}{c(1-c)}(y_k-x^*)+\frac{c+\alpha}{c}(x_{k+1}-x^*)\Big\|^2\\
&=\frac{\alpha}{1-c}\|z_k-x^*\|^2-\frac{\alpha}{c(1-c)}\|y_k-x^*\|^2+\frac{c+\alpha}{c}\|x_{k+1}-x^*\|^2\\
&\quad+\frac{\alpha^2}{c(1-c)^2}\|z_k-y_k\|^2+\frac{\alpha(c+\alpha)}{c^2(1-c)}(\|x_{k+1}-y_k\|^2-c\|x_{k+1}-z_k\|^2)
\end{align*}
and, {using \eqref{e:identz} and \eqref{e:identx}}, we deduce that
\begin{align*}
\|z_{k+1}-x^*\|^2&=\frac{\alpha}{1-c}\|z_k-x^*\|^2-\frac{\alpha}{c(1-c)}\|y_k-x^*\|^2+\frac{c+\alpha}{c}\|x_{k+1}-x^*\|^2\nonumber\\
&\quad +\frac{\alpha^2}{c}\|z_k-x_k\|^2 +\frac{\alpha(c+\alpha)}{c^2}(\|x_{k+1}-x_k\|^2-c\|x_{k}-z_k\|^2)\nonumber\\
&=\frac{\alpha}{1-c}\|z_k-x^*\|^2-\frac{\alpha}{c(1-c)}\|y_k-x^*\|^2+\frac{c+\alpha}{c}\|x_{k+1}-x^*\|^2\nonumber\\
&\quad -\alpha\|z_k-x_k\|^2+\frac{\alpha(c+\alpha)}{c^2}\|x_{k+1}-x_k\|^2,
\end{align*}
which is equivalent to
\begin{align}
\|x_{k+1}-x^*\|^2
&=\frac{c}{c+\alpha}\left(\|z_{k+1}-x^*\|^2-\frac{\alpha}{1-c}\|z_k-x^*\|^2\right)+\frac{c\alpha}{c+\alpha}\|z_k-x_k\|^2\nonumber\\
&\quad+\frac{\alpha}{(c+\alpha)(1-c)}\|y_k-x^*\|^2-\frac{\alpha}{c}\|x_{k+1}-x_k\|^2.
\end{align}
Hence, by replacing in \eqref{e:super} and setting $\lambda_2=1$, we obtain 
\begin{align}
0\ge &\lambda_3 (F(x_{k+1})-F(x^*))- (F(x_{k})-F(x^*))\nonumber\\
&\qquad+\frac{c\lambda_1}{2(c+\alpha)}\left(\rho+\frac{1}{\gamma}\right)\left(\|z_{k+1}-x^*\|^2-\frac{\alpha}{1-c}\|z_k-x^*\|^2\right)\nonumber\\
&\qquad+\left(\frac{\lambda_1}{2}\left(\rho+\frac{1}{\gamma}\right)\frac{c\alpha}{c+\alpha}-\frac{c^2}{2}\left(\frac{1}{\gamma}-\mu\right)\right)\|z_{k}-x_k\|^2\nonumber\\
&\qquad +\frac{\lambda_1}{2}\left(\left(\rho+\frac{1}{\gamma}\right)\frac{\alpha}{(c+\alpha)(1-c)}-\frac{1}{\gamma}+\mu\right)\|y_{k}-x^*\|^2\nonumber\\
&\qquad +\frac{1}{2}\left(\rho+\frac{1}{\gamma}\right)\left(1-\frac{\alpha\lambda_1}{c}\right)
\|x_{k+1}-x_k\|^2\nonumber\\
&\qquad +\frac{\lambda_3}{2}\left(\frac{1}{\gamma}-L\right)\|x_{k+1}-y_k\|^2,
\end{align}
or, equivalently, by defining $\delta_0=\frac{c\lambda_1}{2(c+\alpha)}(\rho+\frac{1}{\gamma})$,
\begin{align}
\label{e:all}
0\ge &\lambda_3\big(F(x_{k+1})-F(x^*)\big)+\delta_0\|z_{k+1}-x^*\|^2\nonumber\\
&\qquad
-\left(F(x_{k})-F(x^*)+\frac{\alpha\delta_0}{1-c}\|z_k-x^*\|^2\right)\nonumber\\
&\qquad+\left(\frac{\lambda_1}{2}\left(\rho+\frac{1}{\gamma}\right)\frac{c\alpha}{c+\alpha}-\frac{c^2}{2}\left(\frac{1}{\gamma}-\mu\right)\right)\|z_{k}-x_k\|^2\nonumber\\
&\qquad +\frac{\lambda_1}{2}\left(\left(\rho+\frac{1}{\gamma}\right)\frac{\alpha}{(c+\alpha)(1-c)}-\frac{1}{\gamma}+\mu\right)\|y_{k}-x^*\|^2\nonumber\\
&\qquad +\frac{1}{2}\left(\rho+\frac{1}{\gamma}\right)\left(1-\frac{\alpha\lambda_1}{c}\right)
\|x_{k+1}-x_k\|^2\nonumber\\
&\qquad +\frac{\lambda_3}{2}\left(\frac{1}{\gamma}-L\right)\|x_{k+1}-y_k\|^2.
\end{align}
Hence, in order to obtain \eqref{e:goal0}, we look for parameters 
satisfying the following inequalities:
\begin{enumerate}
\item\label{e:cond1} $\lambda_3\ge \frac{1-c}{\alpha}$.  
\item\label{e:cond1+} $\kappa:=\frac{\lambda_1}{2}\left(\rho+\frac{1}{\gamma}\right)\frac{c\alpha}{(c+\alpha)}-\frac{c^2}{2}\left(\frac{1}{\gamma}-\mu\right)\ge 0$.  
\item\label{e:cond2} $\left(\rho+\frac{1}{\gamma}\right)\frac{\alpha}{(c+\alpha)(1-c)}-\frac{1}{\gamma}+\mu\ge 0$.
\item\label{e:cond3} $\lambda_1\le \frac{c}{\alpha}$.
\item\label{e:cond4} $\gamma\le \frac{1}{L}$.
\end{enumerate}
{Since, as \eqref{e:all} shows, the linear rate obtained by this 
procedure is $r=1/\lambda_3$, our aim is to choose the parameters 
so as to make $\lambda_3$ as large as possible.}
Note that, since $\lambda_1+1=\lambda_3$, it follows from 
\ref{e:cond1} and \ref{e:cond3} that
\begin{equation}
\frac{\alpha+c}{\alpha}\ge \lambda_1+1=\lambda_3\ge\frac{1-c}{\alpha},
\end{equation}
which {implies} $\alpha\ge 1-2c$ {and 
$\lambda_3\leq\widehat{\lambda}_3:=1+c/\alpha$}. 
{Moreover, for a 
fixed $c$, 
 $\widehat{\lambda}_3$
increases as $\alpha$ decreases, so that the best choice is} 
$\alpha=1-2c$, which implies 
$${\widehat{\lambda}_3=}\lambda_3=\frac{1-c}{\alpha}=\frac{1-c}{1-2c}\quad\text{and}\quad
 \lambda_1=
\frac{c}{\alpha}=\frac{c}{1-2c}.$$
Moreover, condition \ref{e:cond2} becomes
\begin{equation}
\frac{1}{\gamma}\left(\frac{c}{1-c}\right)^2=\frac{1}{\gamma}\left(1-\frac{\alpha}{(c+\alpha)(1-c)}\right)\le \mu+\rho\frac{\alpha}{(c+\alpha)(1-c)}=\mu+\frac{\rho(1-2c)}{(1-c)^2}
\end{equation}
which, combined with \ref{e:cond4}, yields
\begin{equation}
\label{e:interm}
  \frac{1}{L}\ge \gamma\ge \frac{c^2}{\mu(1-c)^2+\rho(1-2c)}.  
\end{equation}
Therefore, by defining $p=\rho/L$ and $q=\mu/L$, we deduce
\begin{equation}
\label{e:corrected}
{1\ge \frac{c^2}{q(1-c)^2+p(1-2c)},}
\end{equation}
{that is, $(1-q)c^2+2(p+q)c-(p+q)\le 0$,}
which implies
\begin{equation}
{c\le \frac{\sqrt{(p+q)(p+1)}-(p+q)}{1-q}
=\frac{\sqrt{p+q}}{\sqrt{p+1}+\sqrt{p+q}}<\frac{1}{2},}
\end{equation}
{the last inequality being equivalent to $q<1$, i.e., to 
$\mu<L$.}
{Since $c\mapsto\lambda_3=(1-c)/(1-2c)$ has derivative 
$(1-2c)^{-2}>0$ on $\left[0,1/2\right[$, the rate $r=1/\lambda_3$ is 
optimized for the largest admissible value of $c$, so we set}
\begin{equation}
\label{e:defc}
{c=\frac{\sqrt{p+q}}{\sqrt{p+1}+\sqrt{p+q}}
=\frac{\sqrt{\mu+\rho}}{\sqrt{L+\rho}+\sqrt{\mu+\rho}}\in\left]0,1\right[.}
\end{equation}
Hence, it follows from 
\eqref{e:interm} that $\gamma=1/L$ and that \ref{e:cond2} is 
satisfied with equality. {Altogether, we deduce from 
$\alpha=1-2c$ and 
$\lambda_3=(1-c)/\alpha$ that}
\begin{equation}
{\alpha=\frac{\sqrt{L+\rho}-\sqrt{\mu+\rho}}
{\sqrt{L+\rho}+\sqrt{\mu+\rho}}
\quad\text{and}\quad
\lambda_3=\frac{\sqrt{L+\rho}}{\sqrt{L+\rho}-\sqrt{\mu+\rho}}
=\frac{1}{r},}
\end{equation}
{with $r$ as in \eqref{e:defr}, and, since 
$\frac{c}{1-c}=\sqrt{\frac{\mu+\rho}{L+\rho}}$, we deduce}
$${\frac{\alpha\delta_0}{1-c}=\frac{c^2(\rho+L)}{2(1-c)^2}
=\frac{\mu+\rho}{2},}$$
\begin{equation}
\label{e:defkappa}
{\kappa=\frac{c^2}{2}\left(\frac{\rho+L}{1-c}-(L-\mu)\right)
=\frac{(\mu+\rho)^{3/2}}
{2\big(\sqrt{\mu+\rho}+\sqrt{L+\rho}\big)}>0,}
\end{equation}
and \eqref{e:all} reduces to
\begin{align}
\label{e:allf}
0&\ge r^{-1} 
\left(F(x_{k+1})-F(x^*)+{\frac{\mu+\rho}{2}}\|z_{k+1}-x^*\|^2\right)\nonumber\\
&\hspace{2cm}
-\left(F(x_{k})-F(x^*)+{\frac{\mu+\rho}{2}}\|z_k-x^*\|^2\right)+\kappa\|z_{k}-x_k\|^2.
\end{align}
The proof is complete.
\end{proof}
\begin{remark}
\label{r:1}
\begin{enumerate}
    \item \label{r:1.1} Note that the function $r=r(\mu,\rho,L)$ in 
    \eqref{e:defr} is strictly decreasing in $\rho$ {(since 
    $\mu<L$)} and 
$r(\mu,0,L)=1-\sqrt{\mu/L}$. Hence, the linear convergence rate is 
strictly smaller when $\rho>0$ and we recover 
\cite[Theorem~4.14]{daspremont2021acceleration} when $\rho=0$ 
and $h=0$.
{In particular, the case $\rho=0$ recovers the parameters of 
V\ref{e:fista}, namely 
$\alpha=(1-\sqrt{\mu/L})/(1+\sqrt{\mu/L})$ and 
$r=1-\sqrt{\mu/L}$ (see, e.g., \cite[Theorem~10.42]{Beck}).}

\item {Observe that the constant $c$ in \eqref{e:defc} is 
strictly positive whenever $\mu+\rho>0$; in particular, 
Theorem~\ref{t:1} covers the case $\mu=0$, in which the strong 
convexity of the objective comes entirely from the non-smooth 
term $h$, and then 
$\alpha=(\sqrt{L+\rho}-\sqrt{\rho})/(\sqrt{L+\rho}+\sqrt{\rho})$ with 
associated rate $r(0,\rho,L)=1-\sqrt{\rho/(L+\rho)}$.}

\item Note that, when $\rho=0$, the constant 
$\kappa$ in \eqref{e:defkappa} reduces to 
$$\kappa=\frac{\mu\sqrt{\mu/L}}{2(1+\sqrt{\mu/L})},$$
which is strictly larger than the constant obtained in the proof of \cite[Theorem~4.14]{daspremont2021acceleration}. This is because we use twice the $\mu$-strong convexity of $f$ in our proof.
\end{enumerate}
\end{remark}

From the Lyapunov-type result in Theorem~\ref{t:1}, we deduce 
the following convergence rates for iterates and values, 
generalizing the result in \cite[Theorem~10.42]{Beck}.
\begin{corollary}
\label{c:rate}
Under the assumptions in Theorem~\ref{t:1}, we have, 
for every $k\in\NN$,
\begin{equation}
\label{e:goal2}
F(x_{k})-F(x^*)\le 
r(\mu,\rho,L)^k\left(F(x_{0})-F(x^*)
+{\frac{\mu+\rho}{2}}\|z_{0}-x^*\|^2\right),
\end{equation}
\begin{equation}
\label{e:goal3}
\|z_{k}-x^*\|\le 
r(\mu,\rho,L)^{k/2}\sqrt{{\frac{2(F(x_{0})-F(x^*))}{\mu+\rho}}
+\|z_{0}-x^*\|^2},
\end{equation}
{and}
\begin{equation}
\label{e:goal4}
{\|x_{k}-x^*\|\le 
r(\mu,\rho,L)^{k/2}\sqrt{\frac{2(F(x_{0})-F(x^*))}{\mu+\rho}
+\|z_{0}-x^*\|^2}.}
\end{equation}
\end{corollary}
\begin{proof}
{The estimates \eqref{e:goal2} and \eqref{e:goal3} follow by 
recurrence from Theorem~\ref{t:1} and \eqref{e:deflyap}, 
together with $\Phi_{\mu,\rho}\ge 0$. For \eqref{e:goal4}, since 
$F$ is $(\mu+\rho)$-strongly convex and $x^*$ minimizes 
$F$, we obtain}
$${\frac{\mu+\rho}{2}\|x_k-x^*\|^2\le F(x_k)-F(x^*)\le 
\Phi_{\mu,\rho}(x_k,z_k)\le 
r(\mu,\rho,L)^k\Phi_{\mu,\rho}(x_0,z_0),}$$
{and the conclusion follows from \eqref{e:deflyap}.}
\end{proof}

\begin{remark}
\label{r:fb}
Suppose that $\alpha=0$ in \eqref{e:fista_z}. Then, for every $k\in\NN$,
$y_k=z_k=x_k$, \eqref{e:fista_z} reduces to the standard 
proximal gradient (or forward-backward) algorithm
\begin{equation}
\label{e:VFB}\tag{FBS}
x_{k+1}=\prox_{\gamma h}\big(x_k-\gamma \nabla f(x_k)\big)
\end{equation}
 and \eqref{e:super} reduces to
\begin{align*}
F(x_{k})-F(x^*)&+\frac{\lambda_1(1-\gamma\mu)}{2\gamma}
\|{x_{k}}-x^*\|^2\\
&\qquad\ge \lambda_3\left(F(x_{k+1})-F(x^*)+\frac{\lambda_1(\rho\gamma+1)}{2\gamma\lambda_3}\|x_{k+1}-x^*\|^2\right)\\
&\qquad\quad 
+\frac{1}{2}\left(\frac{1+\lambda_3}{\gamma}-\lambda_3L+\rho\right)
{\|x_{k+1}-x_k\|^2}.
\end{align*}
By choosing 
\begin{equation}
\lambda_3=\frac{L+\mu+2\rho}{L-\mu}\quad \text{and}\quad
\gamma=\frac{2}{L+\mu},
\end{equation}
and noting that $\lambda_1=\lambda_3-1$, we {obtain that the 
term multiplying $\|x_{k+1}-x_k\|^2$ is zero and we} deduce
\begin{multline}
\label{e:goalb}
F(x_{k+1})-F(x^*)+\frac{\mu+\rho}{2}\|x_{k+1}-x^*\|^2\\\le 
\frac{L-\mu}{L+\mu+2\rho}\left(F(x_{k})-F(x^*)+\frac{\mu+\rho}{2}\|x_{k}-x^*\|^2\right).
\end{multline}
Since the linear convergence rate reduces to 
$\frac{L-\mu}{L+\mu}$ when $\rho=0$, this result complements the 
{results for FBS} obtained in 
\cite{Sigpro1} and {\cite{TayFBS}}. Moreover, {this linear 
convergence rate is better than the one obtained in} 
\cite[Theorem~4.9]{ChPockrev}.
{Observe, incidentally, that the Lyapunov function appearing in 
\eqref{e:goalb} is exactly $\Phi_{\mu,\rho}(x,x)$, with 
$\Phi_{\mu,\rho}$ as in \eqref{e:deflyap}.}
\end{remark}

\section{{Invariance under transfers of strong convexity}}
\label{sec:lev}
For any function $\varphi\colon \HH\to\RX$, we 
denote 
\begin{equation}
(\forall \tau\in\RR)\quad \varphi_{\tau}:=\varphi+\tau\|\cdot\|^2/2.
\end{equation}
Since, for every $\delta\in [-\mu,\rho]$, 
$F=f_{\delta}+h_{-\delta}$, a natural question is if, by choosing a 
specific 
parameter $\delta$, the convergence rates of 
\eqref{e:fista} and/or \eqref{e:FBS} can be improved.
{The results of this section show that this is not the case. 
Indeed, both algorithms with their optimal 
parameters, Lyapunov functions, and linear rates are 
invariant under the transfer $\delta$.}
\begin{proposition}
\label{p:optstrongFB}
Consider the context of Problem~\ref{prob:1}, let 
$\delta\in[-\mu,\rho]$, and set 
\begin{equation}
\label{e:defTdelta}
(\forall \gamma>0)\quad T_{\delta}(\gamma)\colon x\mapsto 
\prox_{\gamma 
h_{-\delta}}(x-\gamma 
\nabla f_{\delta}(x)).
\end{equation}
Then, the following hold:
\begin{enumerate}
\item\label{p:optstrongFBi}
$f_{\delta}$ is differentiable, $(\mu+\delta)$-strongly convex, 
and
$\nabla f_{\delta}$ is $(L+\delta)-$Lipschitz continuous. 

\item\label{p:optstrongFBii} $h_{-\delta}\in\Gamma_0(\HH)$
is $(\rho-\delta)$-strongly convex, for every $\gamma>0$,
$\prox_{\gamma h_{-\delta}}$ is 
$(1+\gamma(\rho-\delta))^{-1}$-Lipschitz continuous, and, { for 
every $\gamma>0$ such that $\gamma\delta<1$,}
$$\prox_{\gamma 
h_{-\delta}}\colon x\mapsto 
\prox_{\frac{\gamma}{1-\gamma\delta}h}
\left(\frac{x}{1-\gamma\delta}\right).$$

\item 
\label{p:optstrongFBiv}
Let $\gamma >0$ and suppose that $\gamma\delta<1$.
Then,
$$T_{\delta}(\gamma)
=T_0\left(\frac{\gamma}{1-\gamma\delta}\right).$$

\end{enumerate}
\end{proposition}
\begin{proof}
\ref{p:optstrongFBi}: {Since $\nabla f_\delta=\nabla f+\delta\Id$, 
the function $f_\delta$ is differentiable. Since 
$\delta\in[-\mu,\rho]$, we have $\mu+\delta\ge0$, 
$L+\delta\ge L-\mu>0$ and}
\begin{equation}
\label{e:stconvmod}
{f_\delta-\frac{\mu+\delta}{2}\|\cdot\|^2
=f-\frac{\mu}{2}\|\cdot\|^2}
\end{equation}
{and }
\begin{equation}
{\frac{L+\delta}{2}\|\cdot\|^2-f_\delta
=\frac{L}{2}\|\cdot\|^2-f.}
\end{equation}
{Since $f$ is 
$\mu$-strongly convex, the right-hand side of \eqref{e:stconvmod} 
is convex and, hence, $f_{\delta}$ is $(\mu+\delta)$-strongly 
convex. Similarly, the $(L+\delta)$-Lipschitz continuous gradient 
is deduced from \cite[Theorem~18.15]{Livre1} and the 
$L$-Lipschitz continuity of $\nabla f$.}

\ref{p:optstrongFBii}: The first part is clear and the strict 
contraction property is obtained from 
\cite[Proposition~23.13]{Livre1}.
Moreover, since $\gamma/(1-\gamma\delta)>0$, it follows from 
\eqref{e:caractprox} and \cite[Example~16.43]{Livre1} that 
\begin{align}
p=\prox_{\gamma h_{-\delta}}x\quad&\Leftrightarrow\quad 
\frac{x-p}{\gamma}\in\partial h_{-\delta}(p)=\partial 
h(p)-\delta p\nonumber\\
&\Leftrightarrow\quad 
\frac{\frac{x}{1-\gamma\delta}-p}{\frac{\gamma}{1-\gamma\delta}}
\in\partial h(p)\nonumber\\
&\Leftrightarrow\quad 
p=\prox_{\frac{\gamma}{1-\gamma\delta}h}
\left(\frac{x}{1-\gamma\delta}\right).
\end{align}

\ref{p:optstrongFBiv}: 
Noting that {$\nabla f_\delta=\nabla f+\delta\Id$} implies that 
$$\Id-\gamma \nabla 
f_{\delta}=(1-\gamma\delta)\left(\Id-\frac{\gamma}{1-\gamma\delta}
\nabla f\right)$$  
and the result follows from \ref{p:optstrongFBii}, $h_0=h$, and 
$f_0=f$.
\end{proof}

{The next result asserts that the optimal linear rate of 
\eqref{e:FBS} 
does not depend on the transfer $\delta$.}
\begin{proposition}
\label{p:contract}
In the context of Proposition~\ref{p:optstrongFB}, 
for every $\gamma\in\left]0,2/(L+\delta)\right[$, 
$T_{\delta}(\gamma)$ is a $\omega(\gamma)$-strict contraction, 
where
\begin{equation}
\omega(\gamma)=\frac{1}{1+\gamma(\rho-\delta)}
\max\{|1-\gamma(\mu+\delta)|,|1-\gamma 
(L+\delta)|\}\in\left]0,1\right[,
\end{equation}
and 
\begin{equation}
\label{e:optpar}
\omega\left(\dfrac{2}{L+\mu+2\delta}\right)
=\min_{\gamma\in[0,2/(L+\delta)]}\omega(\gamma)
=\frac{L-\mu}{L+\mu+2\rho}.
\end{equation}
\end{proposition}
\begin{proof}
Let $\gamma\in\left]0,2/(L+\delta)\right[$. It 
follows from Proposition~\ref{p:optstrongFB}\eqref{p:optstrongFBi} 
and 
\cite[Theorem~2.1]{TayFBS} that
$\Id-\gamma\nabla f_{\delta}$ is a 
$\max\{|1-\gamma(\mu+\delta)|,|1-\gamma(L+\delta)|\}$-strict 
contraction. Hence, we deduce from 
Proposition~\ref{p:optstrongFB}\eqref{p:optstrongFBii} that 
$T_{\delta}(\gamma)=\prox_{\gamma 
h_{-\delta}}\circ (\Id-\gamma\nabla f_{\delta})$ is a 
$\omega(\gamma)$-strict contraction. Noting that 
\begin{equation}
\omega\colon \gamma\mapsto 
\begin{cases}
\dfrac{1-\gamma(\mu+\delta)}{1+\gamma(\rho-\delta)},&\text{if}\:\: 
\gamma\in\left[0,\frac{2}{L+\mu+2\delta}\right[;\\[4mm]
\dfrac{\gamma(L+\delta)-1}{1+\gamma(\rho-\delta)},&\text{if}\:\: 
\gamma\in\left[\frac{2}{L+\mu+2\delta}, \frac{2}{L+\delta}\right]
\end{cases}
\end{equation}
 is continuous, strictly decreasing in $]0,2/(L+\mu+2\delta)[$, and 
 strictly increasing in $]2/(L+\mu+2\delta),2/(L+\delta)[$,
\eqref{e:optpar} follows.
\end{proof}
{In the following result, we conclude that neither the algorithm 
with optimal stepsize nor its linear rate depends 
on $\delta$. Hence, there is no theoretical gain in transferring 
strong convexity from one function to the other in \eqref{e:FBS}.}
\begin{corollary}
\label{c:fbsdelta}
{In the context of Proposition~\ref{p:optstrongFB}, set 
$\gamma_{\delta}=2/(L+\mu+2\delta)$ and let 
$(x_k)_{k\in\NN}$ be generated by 
$x_{k+1}=T_{\delta}(\gamma_{\delta})x_k$, for every $k\in\NN$, 
where $x_0\in\HH$. 
Then $T_{\delta}(\gamma_{\delta})=T_0(2/(L+\mu))$ and}
\begin{equation}
{(\forall k\in\NN)\quad \|x_k-x^*\|\le 
\left(\frac{L-\mu}{L+\mu+2\rho}\right)^{\!\!k}\|x_0-x^*\| .}
\end{equation}
\end{corollary}
\begin{proof}
{Since $\gamma_\delta\delta<1$, 
Proposition~\ref{p:optstrongFB}\eqref{p:optstrongFBiv} gives 
$T_{\delta}(\gamma_\delta)=T_0(\gamma_\delta/
(1-\gamma_\delta\delta))=T_0(2/(L+\mu))$, and the rate follows 
from \eqref{e:optpar} and the Banach--Picard iteration.}
\end{proof}

{We now show that the same invariance holds for 
\eqref{e:fista}. More precisely, transferring strong convexity 
between \(f\) and \(h\) through \(f_{\delta}\) and \(h_{-\delta}\), and 
using the optimal step-size and inertial parameters in 
\eqref{e:param} for the modified problem, leaves both the 
algorithm and its convergence rate unchanged.}
\begin{corollary}
\label{c:2}
In the context of Problem~\ref{prob:1}, let $\delta\in[-\mu,\rho]$,
let $(x_k)_{k\in\NN}$ and $(y_k)_{k\in\NN}$ be the sequences 
generated by the algorithm 
\begin{equation}
\label{e:fistadel}\tag{\ref{e:fista}-$\delta$}
(\forall k\in\NN)\quad 
\begin{cases}
x_{k+1}={\rm prox}_{(L+\delta)^{-1}h_{-\delta}}
\big(y_k-\frac{1}{L+\delta}\nabla f_{\delta}(y_k)\big)\\
y_{k+1}=x_{k+1}+{\frac{\sqrt{(\rho-\delta)+(L+\delta)}
-\sqrt{(\rho-\delta)+(\mu+\delta)}}
{\sqrt{(\rho-\delta)+(L+\delta)}+\sqrt{(\rho-\delta)+(\mu+\delta)}}}
(x_{k+1}-x_k),
\end{cases}
\end{equation}
where $x_0$ and $y_0$ are arbitrary in $\HH$.
Then, {for every $\delta\in[-\mu,\rho]$, \eqref{e:fistadel} 
coincides with \eqref{e:fista} with the parameters \eqref{e:param}. 
In particular, 
for every $k\in\NN$,}
\begin{equation}
{\Phi_{\mu,\rho}(x_{k},z_{k})\le 
\left(1-\sqrt{\frac{\mu+\rho}{L+\rho}}\right)^{\!\! 
k}\Phi_{\mu,\rho}(x_{0},z_{0}),}
\end{equation}
{where $z_k$ is given by \eqref{e:defz} and $\Phi_{\mu,\rho}$ 
is defined in \eqref{e:deflyap}.}
\end{corollary}
\begin{proof}
First note that $f_{\delta}$ and $h_{-\delta}$ satisfy the 
hypotheses of Problem~\ref{prob:1} with 
$\mu_{\delta}=\mu+\delta$, $\rho_{\delta}=\rho-\delta$, and 
$L_{\delta}=L+\delta$, {by 
Proposition~\ref{p:optstrongFB}\eqref{p:optstrongFBi}--%
\eqref{p:optstrongFBii}}. Then, it follows from 
Proposition~\ref{p:optstrongFB}\eqref{p:optstrongFBiv} that 
$T_{\delta}(1/(L+\delta))=T_0(1/L)$ and, therefore, 
the first step of \eqref{e:fistadel} can be written equivalently as
$x_{k+1}=\prox_{h/L}(y_k-\frac{1}{L}\nabla 
f(y_k))$, {which does not depend on $\delta$.
Furthermore, since}
$${\mu_\delta+\rho_\delta=\mu+\rho\quad\text{and}\quad
L_\delta+\rho_\delta=L+\rho,}$$
{and since the constant $c$ in \eqref{e:defc}, the inertial 
parameter $\alpha$ in \eqref{e:param}, the sequence $(z_k)$ in 
\eqref{e:defz}, the Lyapunov function \eqref{e:deflyap}, the 
constant $\kappa$ in \eqref{e:defkappa}, and the rate 
\eqref{e:defr} depend on $(\mu,\rho,L)$ only through $\mu+\rho$ 
and $L+\rho$, all of them are unchanged when $(\mu,\rho,L)$ is 
replaced by $(\mu_\delta,\rho_\delta,L_\delta)$. In particular, the 
inertial step of \eqref{e:fistadel} does not depend on $\delta$ 
either, and \eqref{e:fistadel} coincides with \eqref{e:fista} with the 
parameters \eqref{e:param}. The conclusion now follows from 
Theorem~\ref{t:1} applied to $f_\delta$ and $h_{-\delta}$ and}
\begin{align}
\label{e:goaldel}
{\Phi_{\mu,\rho}(x_{k+1},z_{k+1})}
&{=\Phi_{\mu+\delta,\rho-\delta}(x_{k+1},z_{k+1})}\nonumber\\
&{\le 
r(\mu+\delta,\rho-\delta,L+\delta)
\Phi_{\mu+\delta,\rho-\delta}(x_{k},z_{k})}\nonumber\\
&{=\left(1-\sqrt{\frac{\mu+\rho}{L+\rho}}\right)
\Phi_{\mu,\rho}(x_{k},z_{k})}.
\end{align}
\end{proof}

\begin{remark}
\label{rem:cop}
\begin{enumerate}
\item {Corollaries~\ref{c:fbsdelta} and~\ref{c:2} show that 
neither \eqref{e:FBS} with optimal stepsize nor the inertial scheme 
of Theorem~\ref{t:1} benefits from transferring strong convexity 
between $f$ and $h$. However, the key distinction is whether the 
strong convexity of the non-smooth term is exploited: 
\eqref{e:FBS} 
incorporates it through the step size, whereas \eqref{e:fista} does 
so through the inertial parameter. Indeed, 
V\ref{e:fista}, which ignores $\rho$ and uses 
$\alpha_V=(1-\sqrt{\mu/L})/(1+\sqrt{\mu/L})$, corresponds to 
Theorem~\ref{t:1} with $\rho$ replaced by $0$, and 
Remark~\ref{r:1}\eqref{r:1.1} yields the strictly worse rate 
$r(\mu,0,L)=1-\sqrt{\mu/L}>r(\mu,\rho,L)$ whenever $\rho>0$.}

\item  In order to compare theoretically FBS with the algorithm 
{of Theorem~\ref{t:1}},
note that, without loss of generality, we may assume that 
$L=1$. Indeed, otherwise we can consider the equivalent problem 
with $f_2=f/L$ and $h_2=h/L$, which are $\mu/L$-strongly convex 
with $\nabla f_2$ $1$-Lipschitz continuous and $\rho/L$-strongly 
convex, respectively. By setting, for every $\mu\in[0,1]$ and 
$\rho>0$,
$r_{FBS}(\mu,\rho)=(1-\mu)/(1+\mu+2\rho)$ and 
$r_{\eqref{e:fista}}(\mu,\rho)=1-\sqrt{(\mu+\rho)/(1+\rho)}$, 
we have,  
$${r_{FBS}(\mu,\rho)>r_{\ref{e:fista}}(\mu,\rho)\:\:\Leftrightarrow
\:\:\left(\sqrt{1+\rho}-\sqrt{\mu+\rho}\right)^2>0,}$$ {which is 
always 
true since $\mu<1$.} Hence, 
the 
theoretical upper bound of {the algorithm of 
Theorem~\ref{t:1}}
is always better than that of \eqref{e:FBS}.

\item {On the contrary, if the strong convexity of $h$ is not 
exploited in the inertia, the comparison with \eqref{e:FBS} is no 
longer 
conclusive. In Figure~\ref{fig:comp} we plot 
the regions of $[0,1]\times[0,0.3]$ where 
$r_{FBS}(\mu,\rho)=(1-\mu)/(1+\mu+2\rho)$ is smaller or larger 
than $r_{V\ref{e:fista}}(\mu)=1-\sqrt{\mu}$. We observe that for  
large values of $\rho$, \eqref{e:FBS} has a better upper bound 
than V\ref{e:fista}.}

\begin{figure}
\begin{center}
\includegraphics[scale=0.8]{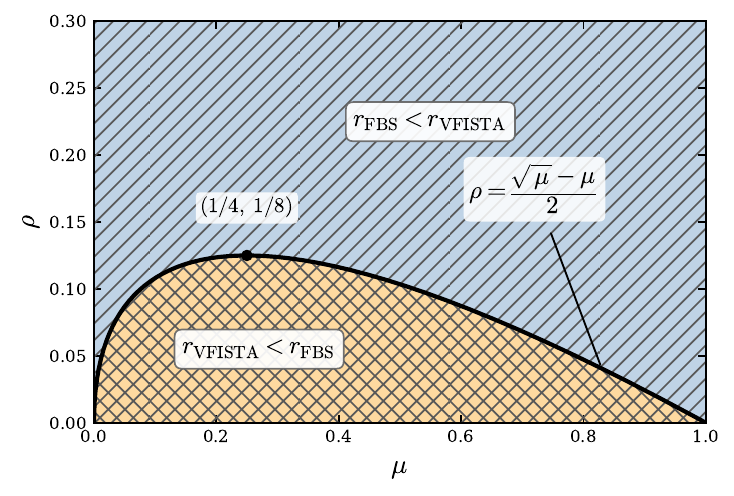}
\end{center}
\caption{Comparison of {V}\ref{e:fista} with FBS depending 
on strong 
convexity parameters $\mu$ and $\rho$.\label{fig:comp}}
\end{figure}

\item Note that {Theorem~\ref{t:1}} provides a simple algorithm 
with fixed inertia and a corresponding autonomous Lyapunov 
function converging linearly at rate $r=1-\sqrt{(\mu+\rho)/(L+\rho)}$.
This rate coincides with that of \cite[Theorem~4.10]{ChPockrev}
for an algorithm with adaptive inertia parameters.
\end{enumerate}

\end{remark}

\section{Numerical experiments}
\label{sec:num}
\begin{figure}[h]
    \centering
    \subfigure[$\mu=1.476\cdot 10^{-6}$]{
\includegraphics[width=0.45\textwidth]
{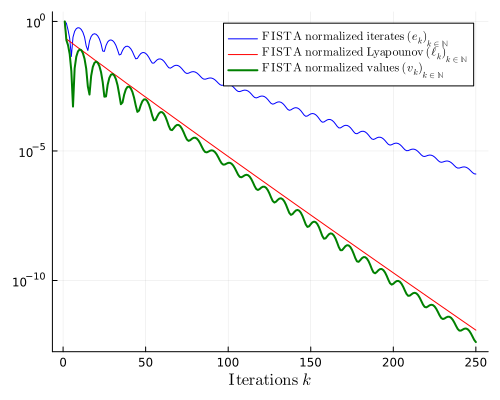}
   }
    \hfill
    \subfigure[$\mu=0.0105$]{
\includegraphics[width=0.45\textwidth]
{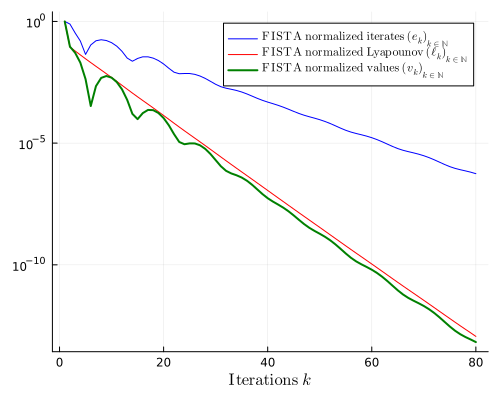}
    
}
    \caption{Normalized error, values, and Lyapunov function for two 
instances with $(n,m)=(50,50)$.\label{sfig:1}}
    \label{fig:1}
\end{figure}

In this section we compare {the algorithm of 
Theorem~\ref{t:1}, the variant V\ref{e:fista} which ignores the 
strong convexity of $h$ in the inertial parameter, and the 
forward-backward 
algorithm \eqref{e:VFB}}
when 
solving strongly convex optimization problems. 
{In view of Corollaries~\ref{c:fbsdelta} and~\ref{c:2}, comparing 
different values of $\delta$ is no longer meaningful, since all the 
resulting algorithms coincide. We therefore compare instead 
inertial parameters built from different amounts of strong 
convexity, namely 
$\alpha_\theta=(\sqrt{L+\theta}-\sqrt{\mu+\theta})/
(\sqrt{L+\theta}+\sqrt{\mu+\theta})$ for 
$\theta\in\{0,\rho/2,\rho\}$, the value $\theta=\rho$ 
corresponding to Theorem~\ref{t:1} and $\theta=0$ to 
V\ref{e:fista}.}
In this context, we 
will consider as an academic example the 
optimization problem
\begin{equation}
\label{e:acadex}
\min_{x\in\RR^n}F(x):=\frac{\rho}{2}\|x+v\|^2+\frac{1}{2}\|Ax-z\|^2, 
\end{equation}
where $(n,m)\in\NN^2$, $A$ is a $m\times n$ real matrix, 
$v\in\RR^n$, $z\in\RR^m$, and $\rho>0$.  
Note that \eqref{e:acadex} 
is a strongly convex optimization problem which can be 
rewritten as
\begin{equation}
\label{e:acadexfh}
\min_{x\in\RR^n}f(x)+h(x), 
\end{equation}
where
\begin{equation}
\begin{cases}
h\colon x\mapsto 
\frac{\rho}{2}\|x+v\|^2\\
f\colon x\mapsto \frac{1}{2}\|Ax-z\|^2.
\end{cases}
\end{equation}
Note that $h$ is $\rho$-strongly convex, $f$ is 
$\lambda_{\min}(A^{\top}A)$-strongly convex, differentiable,
\begin{equation}
\begin{cases}
(\forall \gamma>0)\quad \prox_{\gamma h}\colon x\mapsto 
\dfrac{x-\gamma\rho 
v}{1+\gamma\rho}\\
\nabla f\colon x\mapsto A^{\top}(Ax-z),
\end{cases}
\end{equation}
and $\nabla f$ is $L=\|A^{\top}A\|$-Lipschitz continuous.
Moreover, the optimization problem in \eqref{e:acadex} has the 
unique solution  
\begin{equation}
x^*=(\rho\Id+A^{\top}A)^{-1}(A^{\top}z-\rho v).
\end{equation}
We first fix $(n,m)=(50,50)$, set $\rho=0.1$, 
set $A_0=a\cdot \Id+b\cdot R$, 
$A=A_0/\|A_0^{\top}A_0\|^{1/2}$, 
where
$R$, $v$, and $z$ are generated randomly using the function 
\texttt{rand} in Julia. Since $L=\|A^{\top}A\|=1$, we choose $a$ and 
$b$ in order to 
obtain different values of $\mu$. In particular, we set 
$(a,b)\in\{(0,0.2),(0.58,0.1)\}$, leading to 
$\mu\in\{1.476\cdot 10^{-6}, 0.0105\}$.
In Figure~\ref{fig:1} we plot the normalized error, normalized 
values, and normalized Lyapunov function with respect to 
iterations, which are computed via
\begin{equation}
(\forall k\in\NN)\quad e_k=\frac{\|x_k-x^*\|}{\|x_0-x^*\|},\quad 
v_k=\frac{F(x_k)-F(x^*)}{F(x_0)-F(x^*)},\quad\text{and}\quad
\ell_k=\frac{{\Phi_{\mu,\rho}}(x_k,z_k)}
{{\Phi_{\mu,\rho}}(x_0,z_0)},
\end{equation}
respectively. 

As observed in \cite[Chapter~10.7.4]{Beck} and \cite{suBC},
the error of iterations and value function are not 
monotonically decreasing with respect to iterations. 
On the contrary, Theorem~\ref{t:1} asserts that the sequence 
$\{{\Phi_{\mu,\rho}}(x_k,z_k)\}_{k\in\NN}$ is decreasing, which 
can 
be observed numerically in Figure~\ref{fig:1}. We can interpret 
this as, instead of 
considering either the error in function values or iterates, as it is 
usual, it is preferable to consider $\Phi$ as the
quantity to assess the asymptotic behavior for \eqref{e:fista} 
algorithm, which mixes appropriately errors and value 
functions. 

\begin{figure}[h]
    \centering
    \subfigure[Exact error for $\mu=1.476\cdot 10^{-6}$]{
\includegraphics[width=0.45\textwidth]
{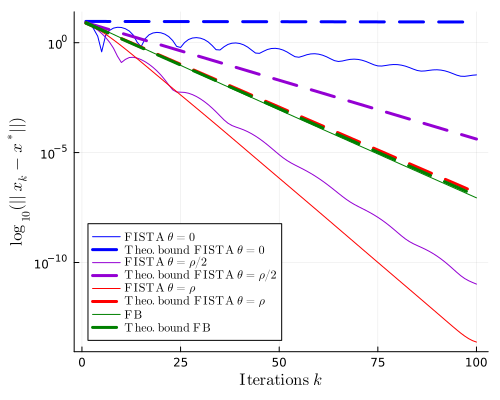}}
    \hfill
    \subfigure[Exact error for $\mu=0.0105$]{
\includegraphics[width=0.45\textwidth]
{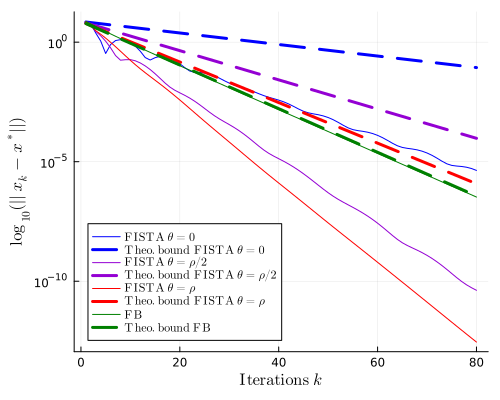}
   }
\\
    \subfigure[Lyapunov function for $\mu=1.476\cdot 10^{-6}$ ]{
\includegraphics[width=0.45\textwidth]
{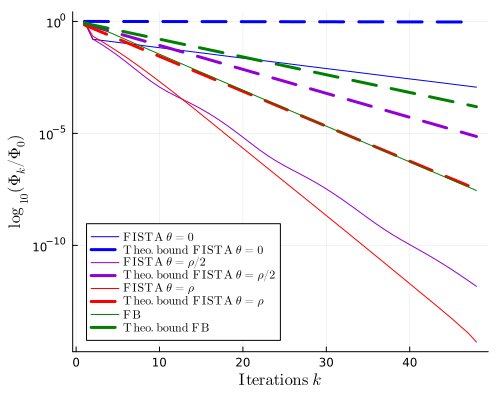}
    
}
    \hfill
    \subfigure[Lyapunov function for $\mu=0.0105$]{
\includegraphics[width=0.45\textwidth]
{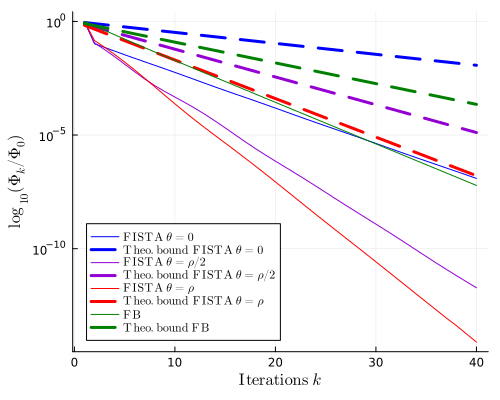}
    
}
    \caption{Exact error and Lyapunov function for FBS and 
    \eqref{e:fista} 
    with $\rho=0.1$ and
{$\theta\in\{0,\rho/2,\rho\}$} for two instances with 
$(n,m)=(50,50)$.\label{sfig:2}}
    \label{fig:2}
\end{figure}

In the same instances, we compare the efficiency of FBS and 
\eqref{e:fista} {when the inertial parameter is built from 
$\theta\in\{0,\rho/2,\rho\}$}. The exact error of 
each algorithm is plotted in Figure~\ref{fig:2} with the 
corresponding theoretical 
upper bounds.
We observe that FBS performs better than {V\ref{e:fista} 
(i.e., $\theta=0$), but when $\theta=\rho$, the algorithm of 
Theorem~\ref{t:1} outperforms} 
all compared algorithms, which is consistent with 
Remark~\ref{rem:cop}.

\begin{figure}[h]
    \centering
    \subfigure[$\mu=1.476\cdot 10^{-6}$]{
\includegraphics[width=0.45\textwidth]
{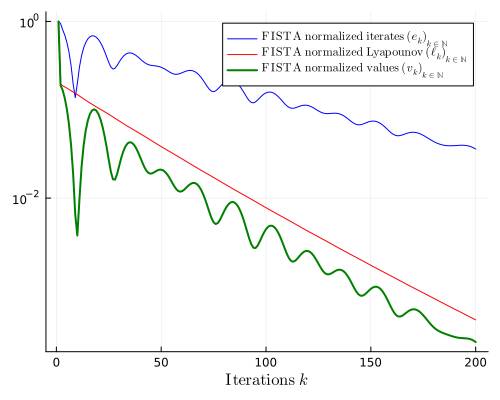}
   }
    \hfill
    \subfigure[$\mu=0.0105$]{
\includegraphics[width=0.45\textwidth]
{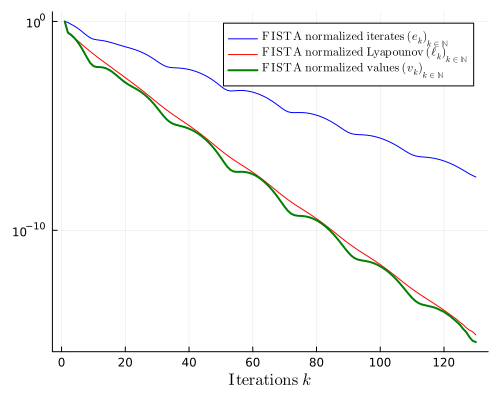}}
\\
    \subfigure[Exact error for $\mu=1.476\cdot 10^{-6}$]{
\includegraphics[width=0.45\textwidth]
{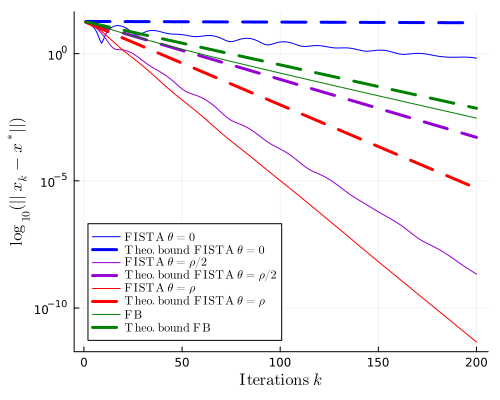}
   }
    \hfill
    \subfigure[Exact error for $\mu=0.0105$]{
\includegraphics[width=0.45\textwidth]
{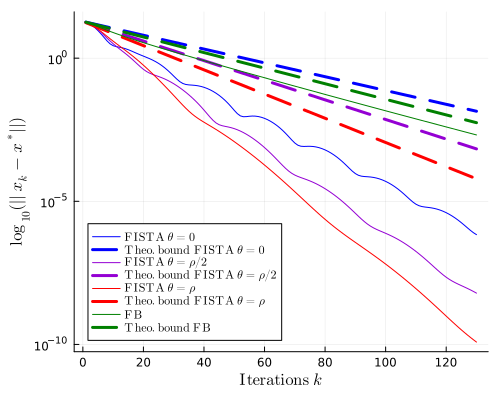}
   }
\\
    \subfigure[Lyapunov function for $\mu=1.476\cdot 10^{-6}$ ]{
\includegraphics[width=0.45\textwidth]
{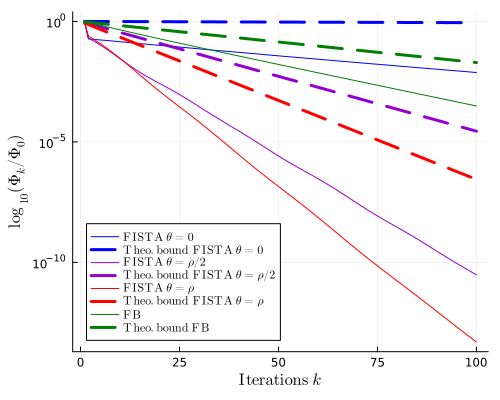}
    
}
    \hfill
    \subfigure[Lyapunov function for $\mu=0.0105$]{
\includegraphics[width=0.45\textwidth]
{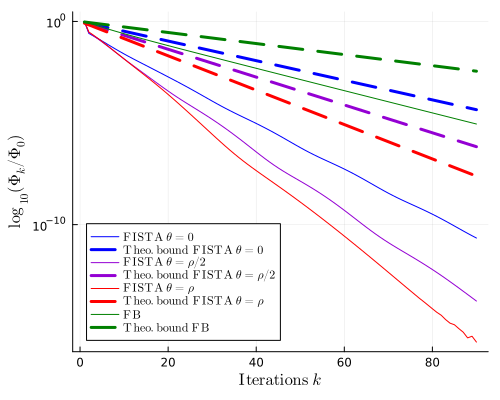}
    
}
    \caption{Exact error and Lyapunov function for FBS and 
    \eqref{e:fista} 
    with $\rho=0.02$,
{$\theta\in\{0,\rho/2,\rho\}$} for two instances with 
$(n,m)=(50,50)$.\label{sfig:3}}
    \label{fig:3}
\end{figure}

In the next experiment, we reduce the strong convexity 
parameter to $\rho=0.02$ for the same randomly generated 
matrices and vectors. In Figure~\ref{fig:3} we observe that 
{\eqref{e:fista} with 
$\theta=0$} now has a better performance in terms of the Lyapunov 
function than FBS for larger $\mu$ but again {the algorithm of 
Theorem~\ref{t:1}} outperforms its competitors.

\begin{figure}[h]
    \centering
    \subfigure[$\mu=2.3\cdot 10^{-9}$]{
\includegraphics[width=0.45\textwidth]
{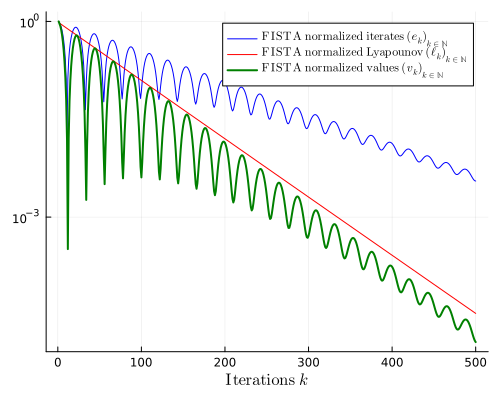}
   }
    \hfill
    \subfigure[$\mu=0.0047$]{
\includegraphics[width=0.45\textwidth]
{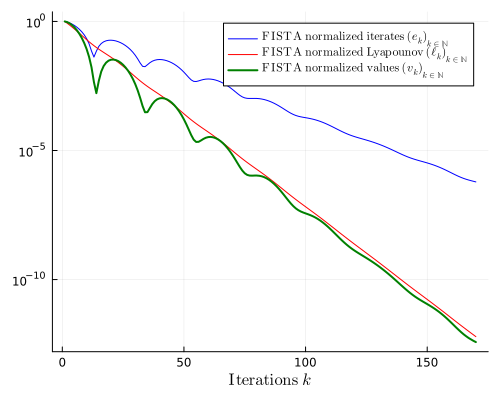}
    
}
\\
    \subfigure[Exact error for $\mu=2.3\cdot 10^{-9}$]{
\includegraphics[width=0.45\textwidth]
{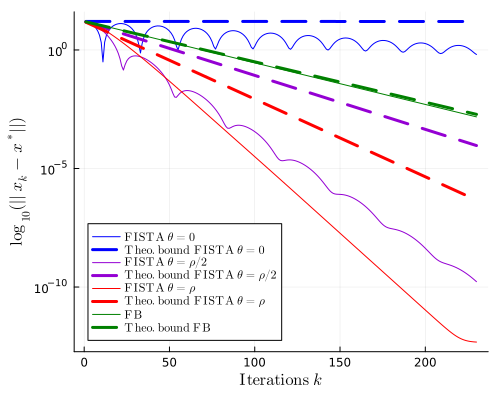}
   }
    \hfill
    \subfigure[Exact error for $\mu=0.0047$]{
\includegraphics[width=0.45\textwidth]
{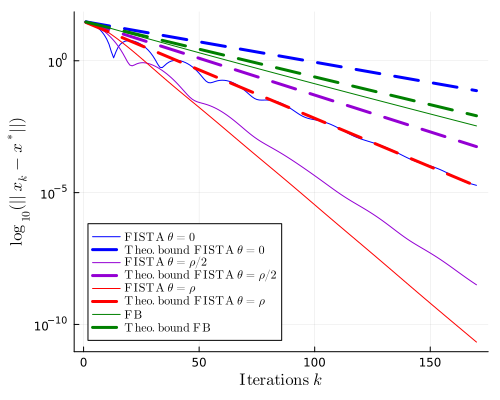}
   }
\\
    \subfigure[Lyapunov function for $\mu=2.3\cdot 10^{-9}$ ]{
\includegraphics[width=0.45\textwidth]
{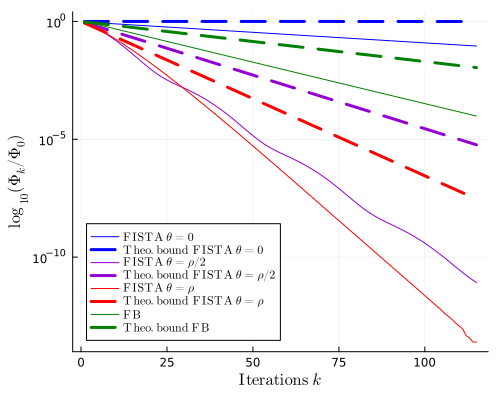}
    
}
    \hfill
    \subfigure[Lyapunov function for $\mu=0.0047$]{
\includegraphics[width=0.45\textwidth]
{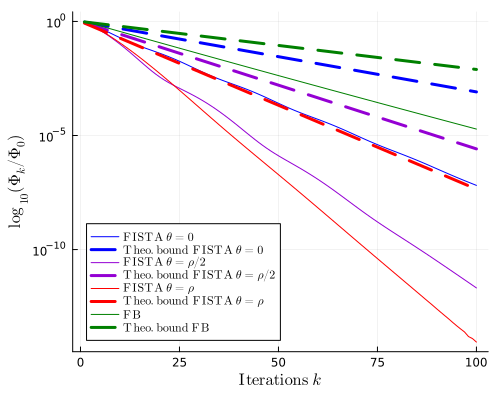}
    
}
    \caption{Exact error and Lyapunov function for FBS and 
    \eqref{e:fista} 
    with $\rho=0.02$,
{$\theta\in\{0,\rho/2,\rho\}$} for two instances with 
$(n,m)=(1000,1000)$.\label{sfig:4}}
    \label{fig:4}
\end{figure}

We finally explore the numerical behavior when we increase the 
dimension to 
$(n,m)=(1000,1000)$ with $(a,b)\in\{(0,0.1), (5,0.1)\}$, leading to 
$\mu\in\{2.3\cdot 10^{-9},0.0047\}$ and keeping 
$\rho=0.02$. In Figure~\ref{fig:4}, we observe pronounced 
oscillations in the values and errors for small values of $\mu$, 
whereas the Lyapunov function decreases monotonically. We also 
observe that 
{the algorithm of Theorem~\ref{t:1}} outperforms its 
competitors and, surprisingly, the 
errors are also monotonically decreasing.

	\bibliographystyle{plain}
	\bibliography{references_2}
\end{document}